\documentclass{amsart}
\usepackage{latexsym, amssymb, amsmath, amsthm, rotating, amscd}
\usepackage{graphics}
\newtheorem{theorem}{Theorem}[section]
\newtheorem{definition}[theorem]{Definition}

\newtheorem{proposition}[theorem]{Proposition}
\newtheorem{corollary}[theorem]{Corollary}
\newtheorem{remark}[theorem]{Remark}

\numberwithin{equation}{section}

\def\Z{\mathbb{Z}}
\def\N{\mathbb{N}}
\def\C{\mathbb{F}}
\def\Q{\mathbb{Q}}

\def\co{\begin{sideways}\begin{sideways}$Y$\end{sideways}\end{sideways}}

\begin{document}

\title[Vertex coalgebras]
{Vertex coalgebras, comodules, cocommutativity and coassociativity}
\author{Keith Hubbard}
\date{\today}

\begin{abstract}
We introduce the notion of vertex coalgebra, a generalization of
vertex operator coalgebras.  Next we investigate forms of
cocommutativity, coassociativity, skew-symmetry, and an
endomorphism, $D^*$, which hold on vertex coalgebras.  The former
two properties require grading. We then discuss comodule structure.
We conclude by discussing instances where graded vertex coalgebras
appear, particularly as related to Primc's vertex Lie algebra and
(universal) enveloping vertex algebras.
\end{abstract}

\maketitle

\section{Introduction}
Since its inception in the 1980s when Borcherds first introduced the
precise notion of a vertex algebra \cite{B}, the theory of vertex
algebras has played an ever expanding role throughout mathematics.
First connected to the understanding of string theory \cite{BPZ},
the Monster group (the largest sporadic finite simple group) and
infinite-dimensional representation theory \cite{FLM}, vertex
algebras now have known connections to modular functions and modular
forms \cite{Z, H3}, Riemann surfaces \cite{H}, Calabi-Yau manifolds
\cite{Fri, GSW}, infinite-dimensional integrable systems
\cite{LW,DJM}, knot and three-manifold invariants \cite{J,RT,W}, and
elliptic cohomology \cite{ST}. Modules of vertex algebras play a
central role in many of these connections.


Vertex operator algebras (VOAs) are a refinement of a vertex
algebras which also encode a conformal structure \cite{FLM}.  (This
structure was actually satisfied by the original example used to
motivate the notion of vertex algebras \cite{B}.)  The dualization
of the VOA structure to that of a vertex operator coalgebra (VOC)
structure has recently been described through the geometry
underlying both structures (\cite{Hub1}, \cite{Hub2}).  In this
paper we generalize the notion of vertex operator coalgebra to the
notion of vertex coalgebra, investigate vertex coalgebra properties,
explore how graded vertex coalgebra structures may be generated from
(and on) graded vertex algebras, and also explore vertex comodule
structures.

Both vertex algebras and VOAs satisfy a natural skew-symmetry as
well as several formal derivative properties. Additionally, the
multiplication operators in these objects satisfy weakened versions
of commutativity (also called locality) and associativity.  These
properties not only prove valuable in understanding the structures
of vertex algebras and VOAs but may also be used to reformulate the
definitions of vertex algebra and VOA in such a way that it is
easier to establish when such a structure exists. Motivated by the
vertex (operator) algebra case, our examination reveals analogous
properties: formal derivatives and skew-symmetry for vertex
coalgebras, and cocommutativity and coassociativity properties for
graded vertex coalgebras. All of these properties extend naturally
to VOCs and in most cases to comodules, allowing us to generate a
number of reformulations of the axioms of vertex coalgebras, graded
vertex coalgebras, VOCs and their comodules.

Classically, one of the most common incarnations of a coalgebra
structure is on the enveloping algebra of a Lie algebra.  We
describe a way to induce this coalgebra structure from the algebra
structure of the enveloping algebra.  We then parallel this
construction for vertex Lie algebras and Primc's enveloping vertex
algebra \cite{P}.  Our construction, which requires imposing a
grading on the vertex Lie algebra, endows the enveloping vertex
algebra with a graded vertex coalgebra structure.  This is a
valuable step in the search for a \emph{general vertex bialgebra},
which subsumes the notion of vertex algebra and vertex coalgebra in
a compatible fashion.

The author thanks his advisor, Katrina Barron, for her insight and
guidance.  He also gratefully acknowledges the financial support of
the Arthur J. Schmidt Foundation.

\section{Definitions and algebraic preliminaries}

We begin, as most papers in vertex algebras, by regurgitating the
required definitions and equations from the calculus of formal
variables.  Both \cite{FLM} and \cite{LL} provide useful expositions
subsuming the discussion below.   All vector spaces will be over a
field $\C$ of characteristic 0.  We will consider $x$, $x_0$, $x_1$
and $x_2$ commuting formal variables throughout, and define the
``formal (or Dirac) $\delta$-function " to be
\begin{equation*}
\delta(x)= \sum_{n \in \Z} x^n.
\end{equation*}
Given $n$ an integer, $(x_1 \pm x_2)^n$ will be understood to be
expanded in nonnegative integral powers of $x_2$.    Note that
$\delta\left(\frac{x_1-x_2}{x_0}\right)$ is a formal power series in
nonnegative powers of $x_2$ (cf. \cite{FHL}, \cite{LL}). The formal
residue, `$Res_{x}$,' refers to the coefficient of the negative
first power of $x$ in the term that follows.  We will use the fact that
\begin{equation} \label{E:Res_x_1}
Res_{x_1} x_1^{-1} \delta \left(\frac{x_2+x_0}{x_1}\right) =1.
\end{equation}

Several basic properties of the $\delta$-function bear mentioning.
First, given a formal Laurent series $X(x_1,x_2) \in (\text{Hom}
(V,W))[[x_1,x_1^{-1},x_2,x_2^{-1}]]$ with coefficients which are
homomorphisms from a vector space $V$ to a vector space $W$, if
$\displaystyle \lim_{x_1 \to x_2} X(x_1,x_2)$ exists (i.e. when
$X(x_1,x_2)$ is applied to any element of $V$, setting $x_1=x_2$
leads to only finite sums in $W$) we have

\begin{equation}\label{E:delta1}
\delta \left( \frac{x_1}{x_2} \right) X(x_1,x_2) = \delta
\left(\frac{x_1}{x_2} \right) X(x_2,x_2).
\end{equation}

\noindent Second, we know that

\begin{equation}\label{E:delta2}
x_1^{-1}\delta \left(\frac{x_2+x_0}{x_1}  \right) =x_2^{-1}\delta
\left(\frac{x_1-x_0}{x_2} \right),
\end{equation}

\noindent and third

\begin{equation}\label{E:delta3}
x_0^{-1}\delta \left(\frac{x_1-x_2}{x_0}  \right) -x_0^{-1}\delta
\left(\frac{x_2-x_1}{-x_0} \right) =x_2^{-1}\delta
\left(\frac{x_1-x_0}{x_2} \right).
\end{equation}

We will also make use of the following equation from Proposition
2.3.6 in \cite{LL}:

\begin{equation} \label{E:split}
\frac{\partial}{\partial x_2} x_2^{-1}\delta \left(\frac{x_1}{x_2}
\right) =(x_1-x_2)^{-2}-(-x_2+x_1)^{-2}.
\end{equation}

%
%

The final `usual suspect' in the onslaught of equations we will not
prove is Taylor's Theorem, which states that  given $V$ a vector
space and $f(x) \in V[[x,x^{-1}]]$,

\begin{equation} \label{E:Taylor}
e^{x_0 \frac{\partial}{\partial x}} f(x)= f(x +x_0),
\end{equation}

\noindent where $e^{x_0 \frac{\partial}{\partial x}}$ is a
formal exponential.

For notational ease, given any vector space $V$, let the linear map $T:V \otimes V \to V \otimes V$
be the transposition operator defined by $T(u \otimes v)=v \otimes u$ for all
$u,v \in V$.

\subsection{The notion of vertex coalgebra} \label{S:vc}
The following description of a vertex coalgebra is the
generalization of a vertex operator coalgebra obtained by omitting
the conformal structure (i.e. the representation of the Virasoro
algebra) on the underlying vector space.


\begin{definition}\label{D:vc}
A \emph{vertex coalgebra} consists of a vector space $V$, together with linear maps

\begin{align*}
\co (x) : V &\mapsto (V \otimes V)[[x,x^{-1}]],
\end{align*}

\begin{equation*}
 c : V \mapsto \C,
\end{equation*}

\noindent called the \emph{coproduct} and the \emph{covacuum map},
respectively, satisfying the following axioms for all $v \in V$:

1. Left Counit:

\begin{equation}\label{E:counit}
(c \otimes Id_V) \co(x)v=v
\end{equation}

2. Cocreation:

\begin{equation} \label{E:cocreat1}
(Id_V \otimes c) \co(x)v \in V[[x]] \ \ \text{and}
\end{equation}
\begin{equation} \label{E:cocreat2}
\lim_{x \to 0} (Id_V \otimes c) \co(x)v=v
\end{equation}

3. Truncation:

\begin{equation} \label{E:trunc}
\co(x)v \in (V \otimes V)((x^{-1}))
\end{equation}

4. Jacobi Identity:

\begin{multline} \label{E:Jac}
x_0^{-1}\delta \left(\frac{x_1-x_2}{x_0} \right) (Id_V \otimes
\co(x_2))  \co(x_1) -x_0^{-1}\delta \left(\frac{x_2-x_1}{-x_0}
\right)
(T \otimes Id_V)\\
 (Id_V \otimes \co(x_1))  \co(x_2)
=x_2^{-1}\delta \left(\frac{x_1-x_0}{x_2} \right) (\co(x_0) \otimes
Id_V)  \co(x_2).
\end{multline}
\end{definition}

The operator $\co(x)$ is linear so that, for example, $(Id_V
\otimes \co(x_2))$ acting on the coefficients of $\co(x_1)v \in (V
\otimes V)[[x_1,x_1^{-1}]]$ is well defined. Notice also, that when
each expression is applied to any element of $V$, the coefficient of
each monomial in the formal variables is a finite sum by Equation (\ref{E:trunc}).
We denote this vertex coalgebra by $V$, or by $(V, \co, c)$ when necessary.

\subsection{Formal derivatives and skew-symmetry on vertex coalgebras} \label{S:Formal_Der}

A reasonable (and fruitful) avenue of investigation is to question the effect
of applying a formal derivative, $d/dx$, to the comultiplication
operator $\co$ of a vertex coalgebra.  The following proposition
gives a convenient description.

\begin{proposition}
Let $(V,\co,c)$ be a vertex coalgebra and $D^*:V \to V$ be the linear
map defined to be

\begin{equation} \label{E:Ddef}
D^*=Res_x x^{-2} (Id_V \otimes c) \co(x).
\end{equation}

\noindent Then

\begin{equation} \label{E:Dco}
(D^* \otimes Id_V)\co(x)=\frac{d}{dx}\co(x).
\end{equation}
\end{proposition}

\begin{proof}
First, we see from the Jacobi identity (\ref{E:Jac}), as well as
(\ref{E:delta2}), that

\begin{multline} \label{E:5171}
Res_{x_0} x_0^{-2} x_1^{-1} \delta \left(\frac{x_2+x_0}{x_1} \right)
(\co(x_0) \otimes Id_V)  \co(x_2) \\
=(x_1-x_2)^{-2} (Id_V \otimes \co(x_2))  \co(x_1) -(x_2-x_1)^{-2}(T
\otimes Id_V) (Id_V \otimes \co(x_1))  \co(x_2).
\end{multline}

\noindent We then use (\ref{E:Ddef}), (\ref{E:Res_x_1}),(\ref{E:5171}),
(\ref{E:counit}), (\ref{E:split}) and (\ref{E:delta1}) in succession
to show that

\begin{align*}
(D^* \otimes &Id_V)\co(x_2)
=Res_{x_0} x_0^{-2} (Id_V \otimes c \otimes Id_V) (\co(x_0) \otimes Id_V)\co(x_2) \\
&=Res_{x_1} Res_{x_0} x_0^{-2} (Id_V \otimes c \otimes Id_V)
x_1^{-1} \delta \left(\frac{x_2+x_0}{x_1}\right) (\co(x_0) \otimes Id_V)\co(x_2) \\
&=Res_{x_1}  (Id_V \otimes c \otimes Id_V) \left((x_1-x_2)^{-2} (Id_V \otimes \co(x_2))  \co(x_1) \right. \\
& \ \ \ \ \ \ \ \ \ \ \left. -(x_2-x_1)^{-2}(T \otimes Id_V) (Id_V \otimes \co(x_1))  \co(x_2) \right) \\
&=Res_{x_1} \left((x_1-x_2)^{-2} \co(x_1)
-(x_2-x_1)^{-2} \co(x_1) \right) \\
&=Res_{x_1} \frac{\partial}{\partial x_2}
x_2^{-1} \delta \left(\frac{x_1}{x_2} \right) \co(x_1) \\
&=Res_{x_1} \frac{\partial}{\partial x_2}
x_2^{-1} \delta \left(\frac{x_1}{x_2} \right) \co(x_2) \\
&=\frac{d}{dx_2}\co(x_2).
\end{align*}
\end{proof}

There is also a convenient way to describe the cocreation axiom in
terms of the operator $D^*$.  First, we note that combining
(\ref{E:counit}) with either (\ref{E:Ddef}) or (\ref{E:Dco})
immediately implies that

\begin{equation} \label{E:cD0}
cD^*=0.
\end{equation}

\noindent Exponentiation of (\ref{E:Dco}), then applying Taylor's Theorem (\ref{E:Taylor})
implies

\begin{align} \label{E:shift}
(e^{x_0 D^*} \otimes Id_V)\co(x)
&=e^{x_0 \frac{\partial}{\partial x}} \co(x) \\
&=\co(x+x_0). \notag
\end{align}

\noindent Composing both sides of (\ref{E:shift}) with the map $Id_V
\otimes c$ gives us a series in only non-negative powers of $x$ so
we set $x=0$ and conclude by (\ref{E:cocreat2}) that

\begin{equation} \label{E:cocreat3}
e^{x_0 D^*} = (Id_V \otimes c) \co(x_0).
\end{equation}

\noindent Equation (\ref{E:cocreat3}) also implies the cocreation
axiom, meaning that the two are equivalent in the presence of the
other axioms.

The operator $D^*$ has a third important application as well.  Vertex
coalgebras possess a natural skew-symmetry via $D^*$ similar to the
case for vertex algebras. More precisely:

\begin{proposition} \label{P:skew-symm}
Given $V$ a vertex coalgebra and $D^*$ as in (\ref{E:Ddef}), we have

\begin{equation} \label{E:skew-symm}
T\co(x)=\co(-x)e^{xD^*}.
\end{equation}
\end{proposition}

\begin{proof}
By two applications of the Jacobi identity (\ref{E:Jac}), we have

\begin{align} \label{E:5201}
&x_2^{-1}\delta \left (\frac{x_1-x_0}{x_2} \right)
(T\co(x_0) \otimes Id_V)  \co(x_2) \\
&=x_0^{-1}\delta \left(\frac{x_1-x_2}{x_0} \right)
(T \otimes Id_V)(Id_V \otimes \co(x_2))  \co(x_1)  \notag \\
& \ \ \ \ \ \ \ \ \ \ \ \ \ \ \ \ \ \ \ \ \ \ \ \ \ \ \ \ \ \ \ \ \
-x_0^{-1}\delta \left(\frac{x_2-x_1}{-x_0} \right)
 (Id_V \otimes \co(x_1))  \co(x_2)  \notag \\
&=x_1^{-1}\delta \left(\frac{x_2+x_0}{x_1} \right) (\co(-x_0)
\otimes Id_V)  \co(x_1). \notag
\end{align}

\noindent But by (\ref{E:delta1}) and (\ref{E:shift}) the right-hand
side of (\ref{E:5201}) is equal to

\begin{multline*}
x_1^{-1}\delta \left(\frac{x_2+x_0}{x_1} \right)
(\co(-x_0) \otimes Id_V)  \co(x_2+x_0) \\
=x_1^{-1}\delta \left(\frac{x_2+x_0}{x_1} \right) (\co(-x_0)e^{x_0
D^*} \otimes Id_V)  \co(x_2).
\end{multline*}

\noindent Taking $Res_{x_1}$ and using (\ref{E:delta2}) we have

\begin{equation*}
(T\co(x_0) \otimes Id_V)  \co(x_2) = (\co(-x_0)e^{x_0 D^*} \otimes
Id_V)  \co(x_2).
\end{equation*}

Finally, apply $(Id_V \otimes Id_V \otimes c)$ to both sides and
note that by cocreation we may set $x_2=0$, getting
(\ref{E:skew-symm}).
\end{proof}

Via skew-symmetry we are now able to interpret the derivative $d/dx$
of $\co$ in a second way, which might be called a \emph{$D^*$-bracket formula for vertex coalgebras}
in view of its similarity to the $D$-bracket formula for vertex algebras.

\begin{proposition} \label{P:Dbracket}
Let $V$ be a vertex algebra and let $D^*$ be as in (\ref{E:Ddef}).
Then

\begin{equation} \label{E:Dbracket}
\frac{d}{dx}\co(x)=\co(x)D^* -(Id_V \otimes D^*) \co(x).
\end{equation}
\end{proposition}

\begin{proof}
Using (\ref{E:skew-symm}), the product rule, (\ref{E:Dco}), and finally
reapplying (\ref{E:skew-symm}),

\begin{align*}
\frac{d}{dx}\co(x) &= \frac{d}{dx} \left( T \co(-x)e^{xD^*} \right) \\
&= T \left(\frac{d}{dx}\co(-x)\right)e^{xD^*} + T \co(-x)\frac{d}{dx}\left(e^{xD^*}\right) \\
&=-T(D^* \otimes Id_V) \co(-x)e^{xD^*} + T \co(-x)e^{xD^*}D^* \\
&=-(Id_V \otimes D^*) \co(x) +  \co(x)D^*.
\end{align*}
\end{proof}

Exponentiating the $D^*$-bracket formula (\ref{E:Dbracket}) yields
\begin{equation} \label{E:Dconjugate1}
(Id_V \otimes e^{-x_0D^*}) \co(x)e^{x_0D^*} = e^{x_0 \frac{d}{dx}}\co(x)
\end{equation}
and applying Taylor's Theorem we have
\begin{equation} \label{E:Dconjugate2}
(Id_V \otimes e^{-x_0D^*}) \co(x)e^{x_0D^*} = \co(x+x_0).
\end{equation}

\subsection{Graded vertex coalgebras}

Although our description of vertex coalgebra parallels Borcherds'
original definition of vertex algebra \cite{B}, some recent
definitions have included a $\Z$-graded vector space (cf.
\cite{F}). Similarly, our focus will be primarily on \emph{graded
vertex coalgebras}.

Let $V = \coprod_{k \in \Z} V_{(k)}$ be a $\Z$-graded vector space
such that $V_{(k)}=0$ for $k<<0$.  Any vector $v \in V_{k}$, $k \in
\Z$, is said to be \emph{homogeneous} and have weight $k$.  We
denote the graded dual space of $V$ by

\begin{equation*}
V'= \coprod_{k \in \Z} V^*_{(k)}=\coprod_{k \in \Z} \text{Hom }(
V_{(k)}, \C),
\end{equation*}

\noindent the algebraic closure of $V$ by

\begin{equation*}
\overline{V}= \prod_{k \in \Z} V_{(k)}=(V')^*,
\end{equation*}

\noindent and the natural pairing of $V'$ with $\overline{V}$ by
$\langle \cdot, \cdot \rangle$. This pairing is also applied to $V$
since it may be viewed as a natural subspace of $\overline{V}$. The
$n$-th tensor product of $V$, denoted $V^{\otimes n}$, is still a
$\Z$-graded vector space (where $v \in V_{k_1} \otimes \ldots
\otimes V_{k_n}$ has weight $k_1 + \ldots + k_n$), and $V^{\otimes
n}$ inherits finite-dimensional homogeneous subspaces if $V$ has
them. Thus $(V^{\otimes n})'$, $\overline{V^{\otimes n}}$ and
$\langle \cdot, \cdot \rangle: (V^{\otimes n})' \times
\overline{V^{\otimes n}} \to \C$ are defined as above.

\begin{definition}
Let $V = \coprod_{k \in \Z} V_{(k)}$ be a $\Z$-graded vector space
such that $V_{(k)}=0$ for $k <<0$.  A vertex coalgebra structure on
$V$ will be said to be \emph{graded} if, given the map $\co (x) : V
\mapsto (V \otimes V)[[x,x^{-1}]]$ and $v \in V_{(r)}$ for some $r
\in \Z$,

\begin{equation*}
\co(x)v = \sum_{k\in \Z} \Delta_k(v) x^{-k-1}
\end{equation*}

\noindent where $\Delta_k(v) \in (V \otimes V)_{(r+k+1)}$ for each
$k\in \Z$. We will say a vertex coalgebra is \emph{finite
dimensionally graded} if it is a graded vertex coalgebra such that
for all $k\in \Z$ we know $\dim V_{(k)}<\infty$.
\end{definition}

In \cite{H}, any Virasoro module that satisfies the condition
$V_{(k)}=0$ for $k$ sufficiently small is called `positive energy
module'.  So what is called a `graded vertex coalgebra' here might
equally well be called a `positive energy graded vertex coalgebra.'
Notice that vertex operator coalgebras are finite dimensionally
graded vertex coalgebras by definition \cite{Hub1}.

Understanding the family of comultiplication operators $\Delta_k:V
\to V \otimes V$ is useful. Already we have a convenient way of
describing the operator $D^*$ of Section \ref{S:Formal_Der} as

\begin{equation*}
D^*=(Id_V \otimes c) \Delta_{-2}.
\end{equation*}

\noindent We may view each operator $\Delta_k$ as a
weight $k+1$ operator in the sense that it maps the $r$-th weight
space of $V$ to the $r+k+1$-st weight space of $V \otimes V$ for
each $r \in \Z$.  Hence $D^*$ is a weight $-1$ operator.

\subsection{Cocommutativity properties} \label{S:cocom}
When investigating non-associative or non-commutative algebras, a
natural question to ask is `How close are they to being associative
or commutative?'. The same question has been studied in detail for
vertex (operator) algebras (cf. \cite{LL}, \cite{FHL}).  In general
the multiplication operator $Y(\cdot,x):V \to (\text{End }V)((x))$
does not commute or associate with itself but does satisfy
conditions called `weak commutativity' (or locality) and `weak
associativity'. To obtain similar results for vertex coalgebras
requires the use of grading and correlations functions. Here we fix
a graded vertex coalgebra $(V, \co, c)$ for consideration.

It might be useful to pause briefly and point out why we do not
consider commutativity and cocommutativity on elements, i.e. $a
\cdot b=b \cdot a$ and $\Delta a = T \Delta a$.  Commutativity or
cocommutativity on elements carries no implication about the
associativity or coassociativity of the elements.  However, if we
generalize slightly and consider (left) commutativity of the
multiplication operator $m$ of an algebra $A$, we see that not only
commutativity of elements, but also associativity is implied.
Specifically, we say that $m$ is left commutative if the diagram

$$
\begin{CD}
A \times A \times A @>Id_A \times m>> A \times A @>m>> A\\
@VT \times Id_A VV  @. @|\\
A \times A \times A @>Id_A \times m>> A \times A @>m>> A.
\end{CD}
$$

\noindent commutes.  This is a stronger requirement than
commutativity of elements, but implies is since $a \cdot b=a \cdot
(b \cdot 1)= b \cdot (a \cdot 1) = b \cdot a$ for all $a,b \in A$.
 It also implies associativity (cf. \cite{LL}). A
comultiplication operator $\Delta$ of a coalgebra $C$ works
similarly: these operators are (right) cocommutative, if the diagram

$$
\begin{CD}
A @>\Delta>> A \times A @>Id_A \times \Delta>> A \times A \times A\\
@| @. @VT \times Id_A VV \\
A   @>\Delta >> A \times A @>Id_A \times \Delta>> A \times A \times
A.
\end{CD}
$$

\noindent commutes.  This property, along with a counit $c$, yields
both cocommutativity:
$$
\Delta= (Id_A \otimes Id_A \otimes c) (Id_A \otimes \Delta) \Delta =
(T \otimes c) (Id_A \otimes \Delta) \Delta = T \Delta,
$$
and coassociativity:
\begin{align*}
(Id_A \otimes \Delta) \Delta
&= (Id_A \otimes T)(Id_A \otimes \Delta) \Delta \\
&=(Id_A \otimes T)(T \otimes Id_A)(Id_A \otimes \Delta) \Delta \\
&=(\Delta \otimes Id_A) T \Delta \\
&=(\Delta \otimes Id_A) \Delta
\end{align*}
where the first and last equalities use cocommutativity of $\Delta$.

Hence we focus our discussion on cocommutativity of comultiplication
operators. Cocommutativity of comultiplication operators does not
hold in all generality for the graded vertex coalgebras, however,
the following weaker result holds. (In Proposition
\ref{P:weakcom_JI}  we will see the extent to which this property
implies coassociativity.)

\begin{proposition} \label{P:weakcom}
\textbf{\emph{(weak cocommutativity)}} Let $v' \in (V^{\otimes
3})'$. There exists $k \in \Z_+$ such that

\begin{equation} \label{E:weakcom}
(x_1-x_2)^k \langle v', (Id_V \otimes \co(x_2))  \co(x_1)v -(T
\otimes Id_V) (Id_V \otimes \co(x_1))  \co(x_2)v \rangle=0
\end{equation}

\noindent for any $v \in V$.
\end{proposition}

\begin{proof}
Since $v' \in (V^{\otimes 3})'$ is finitely generated by elements of
the form $v_1' \otimes v_2' \otimes v_3'$ for $v_1' \in V_{(r)}^*$,
$v_2' \in V_{(s)}^*$, $v_3' \in V_{(t)}^*$, we will only consider
elements of this form. Also pick $N$ such that $V_{(n)}=0$ for all
$n \leq N$.

Multiplying the Jacobi identity by $x_0^k$ and taking $Rex_{x_0}$,
we have

\begin{multline} \label{E:5202}
(x_1-x_2)^k (Id_V \otimes \co(x_2))  \co(x_1)
-(-x_2+x_1)^k (T \otimes Id_V) (Id_V \otimes \co(x_1))  \co(x_2) \\
=Res_{x_0} x_2^{-1}\delta \left(\frac{x_1-x_0}{x_2} \right) x_0^k
(\co(x_0) \otimes Id_V)  \co(x_2).
\end{multline}

\noindent We know that for some $M \in \Z_+$

\begin{equation} \label{E:5204}
\langle v_1' \otimes v_2', \co(x_0)v \rangle \in \C x_0^{-M}[[x_0]]
\end{equation}

\noindent for any choice of $v \in V_{(q)}$, since a non-zero
coefficient implies that $\text{wt } v_1' + \text{wt } v_2' =
\text{wt }\Delta_n(v) = \text{wt } v + n +1$ or $ r+s=q+n+1$, which
implies that $r+s>N+n+1$ or $-n-1>N-r-s$. Simply let $M \geq -N+r+s$.

Because the choice of $M$ is independent of the vector $\co(x_0)$
acts on, for all $k>M$ we have

\begin{equation*}
Res_{x_0} x_2^{-1}\delta \left(\frac{x_1-x_0}{x_2} \right) x_0^k
\langle v_1' \otimes v_2' \otimes v_3', (\co(x_0) \otimes Id_V)
\co(x_2)v \rangle=0
\end{equation*}

\noindent for any choice of $v \in V$.  Hence (\ref{E:5202}) implies
the proposition.
\end{proof}

We can understand a great deal about compositions of $\co$ operators
by investigating how to obtain them from rational functions.  Our
discussion parallels the discourse in \cite{LL} on expansions into
formal Laurent series. Let $\C[x_1,x_2]_S$ be the ring of rational
functions obtained by inverting the products of (zero or more)
elements of $S$, where $S$ is the set of nonzero homogeneous linear
polynomials in $x_1$ and $x_2$.  Let $\iota_{1 2}:\C[x_1,x_2]_S \to
\C[[x_1,x_1^{-1},x_2,x_2^{-1}]]$ be defined by mapping an element
$\frac{g(x_1,x_2)}{x_1^r x_2^s (x_1-x_2)^t}$ to
$\frac{g(x_1,x_2)}{x_1^r x_2^s}$ times $\frac{1}{(x_1-x_2)^t}$
expanded in nonnegative powers of $x_2$ where $g(x_1,x_2) \in
\C[x_1,x_2]$ and $r,s,t \in \N$.  The map $\iota_{1 2}$ is injective
and hence can be inverted on $\iota_{1 2}[\C[x_1, x_2]_S].$ This
generalizes to finitely many formal variables, i.e. $\iota_{1 \cdots
n}$, and any ordering of the formal variables (cf. Section 3.1 of
\cite{FHL}).

\begin{proposition}\label{P:rr&com} (a) \textbf{\emph{(right rationality)}}
Let $v' \in (V^{\otimes 3} )'$ and $v \in V$. Then the formal series
$\langle v', (Id_V \otimes \co(x_2)) \co(x_1)v \rangle$ is in
$\C[x_1 ,x_2^{-1}][[x_1^{-1},x_2]]$ and in fact

\begin{equation} \label{E:5211}
\langle v', (Id_V \otimes \co(x_2)) \co(x_1)v \rangle = \iota_{1 2}
f(x_1,x_2)
\end{equation}

\noindent where $f(x_1,x_2) \in \C[x_1,x_2]_S$ is uniquely
determined and is of the form

\begin{equation} \label{E:129}
f(x_1,x_2)= \frac{g(x_1,x_2)}{x_1^r x_2^s (x_1-x_2)^t}
\end{equation}

\noindent for some $g(x_1,x_2) \in \C[x_1,x_2]$ and $r,s,t \in \N$.

(b) \textbf{\emph{(cocommutativity)}} It is also the case that

\begin{equation*}
\iota^{-1}_{1 2} \langle v', (Id_V \otimes \co(x_2)) \co(x_1)v
\rangle =\iota^{-1}_{2 1} \langle v', (T \otimes Id_V) (Id_V \otimes
\co(x_1)) \co(x_2)v \rangle.
\end{equation*}
\end{proposition}

\begin{proof}
Let $t\in \Z_+$ be chosen relative to $v'$ so that weak
cocommutativity holds, i.e.

\begin{multline} \label{E:5203}
(x_1-x_2)^t \langle v',  (Id_V \otimes \co(x_2)) \co(x_1)v \rangle \\
=(x_1-x_2)^t \langle v', (T \otimes Id_V)
 (Id_V \otimes \co(x_1))  \co(x_2)  v \rangle.
\end{multline}

\noindent This is an equality of formal power series in
$\C[[x_1,x_1^{-1},x_2,x_2^{-1}]]$, but the left-hand side has only
finitely many positive powers of $x_1$ by truncation (\ref{E:trunc})
and only finitely many negative powers of $x_2$ by (\ref{E:5204}).
Similarly the right-hand side of (\ref{E:5203}) has only finitely
many positive powers of $x_2$ and finitely many negative powers of
$x_1$. Hence both sides of (\ref{E:5203}) must be equal to a Laurent
polynomial (unique for a given $t$)

\begin{equation*}
\frac{g(x_1,x_2)}{x_1^r x_2^s}
\end{equation*}

\noindent for some $g(x_1,x_2) \in \C[x_1,x_2]$, with $r,s \in \Z$.
Let $f(x_1,x_2)$ be the rational function uniquely defined by

\begin{equation*}
f(x_1,x_2) =\frac{g(x_1,x_2)}{x_1^r x_2^s (x_1-x_2)^t}.
\end{equation*}

\noindent Since neither $(x_1-x_2)^t$ nor $\langle v', (Id_V \otimes
\co(x_2)) \co(x_1)v \rangle$ have infinitely many negative powers of
$x_2$, we may multiply by $(x_1-x_2)^{-t} \in \C[[x_1^{-1},x_2]]$ so
that

\begin{align*}
\langle v',  (Id_V \otimes \co(x_2)) \co(x_1)v \rangle
&=(x_1-x_2)^{-t} (x_1-x_2)^t \langle v',  (Id_V \otimes \co(x_2)) \co(x_1)v \rangle \\
&= \iota_{1 2} f(x_1,x_2).
\end{align*}

\noindent (See Section 2.1 on \cite{LL} for a discussion of when
such products are defined.) Thus part (a) is satisfied.  Similarly,
using finite negative powers of $x_1$ of the right-hand side of (\ref{E:5203}),
 we observe that

\begin{equation*}
\langle v',  (T \otimes Id_V)(Id_V \otimes \co(x_1)) \co(x_2)v
\rangle = \iota_{2 1} f(x_1,x_2)
\end{equation*}

\noindent satisfying part (b).
\end{proof}

\begin{remark} \label{R:wcom_com}
The proof of Proposition \ref{P:rr&com} actually demonstrates that
right rationality and cocommutativity follow from the axioms of a
graded vertex coalgebra without the Jacobi identity, but with the
addition of weak cocommutativity (\ref{E:weakcom}).
\end{remark}

\subsection{Coassociativity properties}
In the same way as we have investigated cocommutativity, we now
investigate coassociativity and achieve the following results.

\begin{proposition} \label{P:weakassoc}
\textbf{\emph{(weak coassociativity)}} Let $v' \in (V^{\otimes
3})'$. There exists $k \in \Z_+$ such that

\begin{equation} \label{E:weakassoc}
(x_0+x_2)^k \langle v', (\co(x_0) \otimes Id_V)  \co(x_2)v -(Id_V
\otimes \co(x_2))  \co(x_0 + x_2)v \rangle=0
\end{equation}

\noindent for any $v \in V$.
\end{proposition}

\begin{proof}
Since $v' \in (V^{\otimes 3})'$ is finitely generated by elements of
the form $v_1' \otimes v_2' \otimes v_3'$ for $v_1' \in V_{(r)}^*$,
$v_2' \in V_{(s)}^*$, $v_3' \in V_{(t)}^*$, we again consider
elements of this form and again pick $N$ such that $V_{(n)}=0$ for
all $n \leq N$.

To achieve weak cocommutativity we took an appropriate residue of
$x_0$ from the Jacobi identity and noted that the right-hand side
had only finitely many negative $x_0$ terms. We will now use the
same approach with $x_1$ and the second term of the Jacobi identity.
Multiplying the Jacobi identity by $x_1^k$ and taking $Rex_{x_1}$,
we have

\begin{multline} \label{E:5212}
(x_0+x_2)^k (Id_V \otimes \co(x_2))  \co(x_0+x_2)
-Res_{x_1} x_0^{-1}\delta \left(\frac{x_2-x_1}{-x_0} \right) \\
x_1^k (T \otimes Id_V) (Id_V \otimes \co(x_1))  \co(x_2) =(x_2
+x_0)^k (\co(x_0) \otimes Id_V) \co(x_2)
\end{multline}

\noindent by applying $(\ref{E:delta2})$ and $(\ref{E:delta1})$ to
the first term, and $(\ref{E:delta2})$ to the last term.  We know
that for some $M \in \Z_+$

\begin{equation} \label{E:5213}
\langle v_1' \otimes v_3', \co(x_1)v \rangle \in \C x_1^{-M}[[x_1]]
\end{equation}

\noindent for any choice of $v \in V$, as in (\ref{E:5204}). Thus
for $k>M$ we see that

\begin{equation*}
Res_{x_1} x_0^{-1}\delta \left(\frac{x_2-x_1}{-x_0} \right) x_1^k
\langle v_1' \otimes v_2' \otimes v_3', (T \otimes Id_V) (Id_V
\otimes \co(x_1))  \co(x_2) v \rangle=0
\end{equation*}

\noindent for any choice of $v \in V$, so (\ref{E:5212}) implies the
proposition.
\end{proof}

Again mirroring our above discussion, we investigate an
interpretation of coassociativity in terms of rational functions.

\begin{proposition}\label{P:lr&ass} (a) \textbf{\emph{(left rationality)}}
Let $v' \in (V^{\otimes 3})'$ and $v \in V$. Then the formal series
$\langle v', (\co(x_0) \otimes Id_V) \co(x_2)v \rangle$ is in
$\C[x_0^{-1},x_2][[x_0,x_2^{-1}]]$ and in fact

\begin{equation*}
\langle v', (\co(x_0) \otimes Id_V) \co(x_2)v \rangle = \iota_{2 0}
k(x_0,x_2)
\end{equation*}

\noindent where $k(x_0,x_2) \in \C[x_0,x_2]_S$ is uniquely
determined and is of the form

\begin{equation*}
k(x_0,x_2)= \frac{h(x_0,x_2)}{x_0^r x_2^s (x_0+x_2)^t}
\end{equation*}

\noindent for some $h(x_0,x_2) \in \C[x_0,x_2]$ and $r,s,t \in \Z$.

(b) \textbf{\emph{(coassociativity)}} It is also the case that
\begin{equation*}
\iota^{-1}_{1 2} \langle v', (Id_V \otimes \co(x_2)) \co(x_1)v
\rangle =\left. \left( \iota^{-1}_{2 0} \langle v', (\co(x_0)
\otimes Id_V) \co(x_2)v \rangle \right) \right|_{x_0=x_1-x_2}.
\end{equation*}

\end{proposition}

\begin{proof}
Let $t\in \Z_+$ be chosen relative to $v'$ so that weak coassociativity
holds, i.e.

\begin{multline} \label{E:5214}
(x_0+x_2)^t \langle v', (\co(x_0) \otimes Id_V)  \co(x_2)v \rangle
=(x_0+x_2)^t \langle v',(Id_V \otimes \co(x_2))  \co(x_0 + x_2)v
\rangle.
\end{multline}

\noindent The left-hand side of (\ref{E:5214}) has only finitely
many positive powers of $x_2$ by truncation (\ref{E:trunc}) and only
finitely many negative powers of $x_0$ by (\ref{E:5213}). The right-hand
side of (\ref{E:5214}), on the other hand, has only finitely
many positive powers of $x_0$ and finitely many negative powers of
$x_2$.  Hence both sides of (\ref{E:5214}) must be equal to a
Laurent polynomial (unique given $t$)

\begin{equation*}
\frac{h(x_0,x_2)}{x_0^r x_2^s}
\end{equation*}

\noindent for some $h(x_0,x_2) \in \C[x_0,x_2]$, with $r,s \in \Z$.
Let $k(x_0,x_2)$ be the rational function defined by

\begin{equation*}
k(x_0,x_2) =\frac{h(x_0,x_2)}{x_0^r x_2^s (x_0+x_2)^t}.
\end{equation*}

\noindent Since neither $(x_0+x_2)^t$ nor $\langle v', (\co(x_0)
\otimes Id_V)  \co(x_2)v \rangle$ have infinitely many negative
powers of $x_0$, we may multiply by
$(x_2+x_0)^{-t}\in\C[[x_2^{-1},x_0]]$ so that

\begin{align*}
\langle v', (\co(x_0) \otimes Id_V)  \co(x_2)v \rangle
&=(x_2+x_0)^{-t} (x_0+x_2)^t \langle v', (\co(x_0) \otimes Id_V)  \co(x_2)v \rangle \\
&= \iota_{2 0} k(x_0,x_2).
\end{align*}

\noindent This satisfies part (a).  Since the right-hand side of
(\ref{E:5214}) has only finite negative powers of $x_2$, we may
similarly verify that

\begin{align*}
\langle v', (Id_V \otimes \co(x_2))  \co(x_0 + x_2)v \rangle =
\iota_{0 2} k(x_0,x_2).
\end{align*}

\noindent Thus

\begin{align*}
\iota^{-1}_{1 2} \langle v', (Id_V \otimes \co(x_2)) \co(x_1)v
\rangle &=\left. \iota^{-1}_{0 2}
      \langle v', (Id_V \otimes \co(x_2))  \co(x_0 + x_2)v \rangle\right|_{x_0=x_1-x_2}\\
&=\left. k(x_0,x_2)\right|_{x_0=x_1-x_2}\\
&=\left. \left( \iota^{-1}_{2 0} \langle v', (\co(x_0) \otimes Id_V)
\co(x_2)v \rangle \right) \right|_{x_0=x_1-x_2}.
\end{align*}
\end{proof}

\begin{remark} \label{R:wassoc_assoc}
The proof of Proposition \ref{P:lr&ass} actually demonstrates that
left rationality and coassociativity follow from the axioms of a
graded vertex coalgebra without the Jacobi identity, but with the
addition of weak coassociativity (\ref{E:weakassoc}).
\end{remark}

\section{Equivalent characterizations of graded vertex coalgebras}
As we saw at the beginning of Section \ref{S:cocom}, in classical
algebra cocommutativity of comultiplication operators implies
cocommutativity \emph{and} coassociativity.  An analog of this
result is true for graded vertex coalgebras.  We use these
observations to investigate a way to replace the Jacobi identity in
the definition of graded vertex coalgebra.

Although, there is good reason to keep the Jacobi identity in the
definition of vertex algebra or coalgebra (for instance that of
keeping the definitions of a vertex (co)algebra and a (co)module
over a vertex (co)algebra in the same form), it is useful to have
equivalent characterizations of the Jacobi identity. In this section
we investigate two equivalent characterizations of the Jacobi
identity, and hence alternate characterizations of graded vertex
coalgebras in general. Note that these characterizations do not
apply to vertex coalgebras in general because (\ref{E:5204}) and
(\ref{E:5213}) arise from the grading.

\subsection{The Jacobi identity from cocommutativity and the $D^*$-bracket}

We begin by claiming Proposition 7.4.1 from \cite{Hub1}. The
statement there is about the Jacobi identity being replaced in the
context of a vertex operator coalgebra, but the proof is void of
any reference to conformal structure (i.e. a representation of the
Virasoro algebra) and is applicable to graded vertex coalgebras
without modification.

\begin{proposition}\label{P:RCAtoJI}
The Jacobi identity may be replaced in the definition of a graded
vertex algebra by right and left rationality, cocommutativity and
coassociativity, i.e. the claims of Propositions \ref{P:rr&com} and
\ref{P:lr&ass}.
\end{proposition}

Putting together several of the preceding results, it is now
possible to argue that the Jacobi identity is a consequence of weak
cocommutativity and the $D^*$-bracket.

\begin{proposition} \label{P:weakcom_JI}
In the presence of all the axioms of a graded vertex coalgebra except
the Jacobi identity, weak cocommutativity (\ref{E:weakcom}) and the
$D^*$-bracket (\ref{E:Dbracket}) imply the Jacobi identity.
\end{proposition}

\begin{proof}
First, we will show that Equation (\ref{E:cocreat3}) and
skew-symmetry (\ref{E:skew-symm}) hold without use of the Jacobi identity, and then show
that weak coassociativity also holds.  Via Remarks
\ref{R:wcom_com} and \ref{R:wassoc_assoc}, that means that right and
left rationality, cocommutativity and coassociativity all hold.
Thus, by Proposition \ref{P:RCAtoJI}, the Jacobi identity holds.

Since the $D^*$-bracket immediately implies (\ref{E:Dconjugate2}), by
(\ref{E:cD0}) we have

\begin{align*}
(Id_V \otimes c)\co(x + x_0)
&=(Id_V \otimes c e^{-x_0 D^*}) \co(x) e^{x_0 D^*} \\
&=(Id_V \otimes c) \co(x) e^{x_0 D^*}.
\end{align*}

\noindent By the cocreation axiom, we may set $x=0$ and get
(\ref{E:cocreat3}), i.e.

\begin{equation} \label{E:newcocreat3}
(Id_V \otimes c)\co(x_0) = e^{x_0 D^*}.
\end{equation}

\noindent Now let $v' \in (V^{\otimes 2})'$ and choose $k\in \Z_+$
so that weak cocommutativity holds on $v' \otimes c$. Applying
(\ref{E:newcocreat3}), then (\ref{E:Dconjugate2}) to the right-hand
side of weak commutativity, we have

\begin{align*}
(x_1-x_2)^k \langle v' \otimes c, & (Id_V \otimes
\co(x_2))\co(x_1)v \rangle \\
&= (x_1-x_2)^k \langle v' \otimes c, (T \otimes Id_V)(Id_V \otimes \co(x_1))\co(x_2)v \rangle \\
&= (x_1-x_2)^k \langle v' , T (Id_V \otimes e^{x_1 D^*}) \co(x_2)v \rangle \\
&= (x_1-x_2)^k \langle v' , T \co(x_2-x_1)e^{x_1 D^*} v \rangle.
\end{align*}

\noindent But since $(Id_V \otimes c)\co(x_2) \in V[[x_2]]$, we may
take $x_2=0$ and get

\begin{equation*}
x_1^k \langle v', \co(x_1)v \rangle = x_1^k \langle v' , T
\co(-x_1)e^{x_1 D^*} v \rangle.
\end{equation*}

\noindent Multiplying by $x_1^{-k}$ is well-defined so we get

\begin{equation} \label{E:5215}
\langle v', \co(x_1)v \rangle = \langle v' , T \co(-x_1)e^{x_1 D^*} v
\rangle,
\end{equation}

\noindent a form of skew-symmetry (independent of $v'$ and $v$).

Finally, let $v' \in (V^{\otimes 3})'$ and choose $k$ such that weak
commutativity holds for $v'(Id_V \otimes T)$.  Employing
(\ref{E:5215}), (\ref{E:Dconjugate2}), weak cocommutativity, and
(\ref{E:5215}), in order, we have

\begin{align*}
(x_0+x_2)^k &\langle v', (Id_V \otimes \co(x_2)) \co(x_0+x_2)v \rangle \\
&=(x_0+x_2)^k \langle v', (Id_V \otimes T \co(-x_2)e^{x_2 D^*}) \co(x_0+x_2)v \rangle \\
&=(x_0+x_2)^k \langle v', (Id_V \otimes T \co(-x_2)) \co(x_0) e^{x_2 D^*}v \rangle \\
&=(x_0+x_2)^k \langle v', (Id_V \otimes T) (T \otimes Id_V)(Id_V \otimes \co(x_0)) \co(-x_2) e^{x_2 D^*}v \rangle \\
&=(x_0+x_2)^k \langle v', (Id_V \otimes T) (T \otimes Id_V)(Id_V \otimes \co(x_0)) T\co(x_2)v \rangle \\
&=(x_0+x_2)^k \langle v', (\co(x_0) \otimes Id_V)  \co(x_2)  v
\rangle.
\end{align*}

\noindent This is weak coassociativity, completing our proof.
\end{proof}

Proposition \ref{P:weakcom_JI} allows us to give the following
equivalent definition of a graded vertex coalgebra.

\begin{definition} \label{D:gvc_alt}
A \emph{graded vertex coalgebra} may be described as a collection of
data

\noindent
\begin{itemize}
\item a $\Z$-graded vector space $V=\coprod_{k = -N}^{\infty} V_{(k)}$, with $N \in \N$
\item a linear map $c : V \mapsto \C$
\item a linear map

\begin{align*}
\co (x) : V &\to (V \otimes V)((x^{-1})) \\
v &\mapsto \co(x)v = \sum_{k\in \Z} \Delta_k(v) x^{-k-1},
\end{align*}
such that for $v \in V_{(j)}$, $\Delta_k(v) \in (V \otimes
V)_{(j+k+1)}$,
\end{itemize}

\noindent
satisfying the following axioms for $v \in V$ and $v' \in (V^{\otimes 3})'$: \\
1. $(c \otimes Id_V) \co(x)=Id_V$ \\
2. $(Id_V \otimes c) \co(x)v \in V[[x]]$, and

$\lim_{x \to 0} (Id_V \otimes c) \co(x)v=v$ \\
3. $\frac{d}{dx}\co(x)=\co(x)D^* -(Id_V \otimes D^*) \co(x)$ where $D^*=(Id_V \otimes c)\Delta_{-2}$ \\
4. There exists $k \in \Z_+$ depending on $v'$ but independent of $v$ such that

\begin{equation*}
(x_1-x_2)^k \langle v', (Id_V \otimes \co(x_2))  \co(x_1)v -(T
\otimes Id_V) (Id_V \otimes \co(x_1))  \co(x_2)v \rangle=0.
\end{equation*}
\end{definition}
Here, we combine the truncation axiom into the definition of the
range of $\co(x)$.

\subsection{The Jacobi identity from coassociativity and skew-symmetry}
There is another equivalent substitute for the Jacobi identity
besides weak cocommutativity and the $D^*$-bracket, which also
yields an equivalent characterization of graded vertex coalgebra.

\begin{proposition} \label{P:weakassoc_JI}
In the presence of all the axioms of a graded vertex coalgebra except
the Jacobi identity, weak coassociativity (\ref{E:weakassoc}) and
skew-symmetry (\ref{E:skew-symm}) imply the Jacobi identity.
\end{proposition}

\begin{proof}
Mirroring the proof of Proposition \ref{P:weakcom_JI}, we need only
prove that weak cocommutativity is satisfied.  First, note that the
left counit axiom (\ref{E:counit}) and skew-symmetry
(\ref{E:skew-symm}) imply (\ref{E:cocreat3}). Now let $v' \in
(V^{\otimes 2})'$ and choose $k\in \Z_+$ so that weak
coassociativity holds on $v' \otimes c$:

\begin{equation*}
(x_0+x_2)^k \langle v' \otimes c, (\co(x_0) \otimes Id_V)  \co(x_2)v
\rangle =(x_0+x_2)^k \langle v' \otimes c, (Id_V \otimes \co(x_2))
\co(x_0 + x_2)v \rangle.
\end{equation*}

\noindent Adding (\ref{E:cocreat3}), we get

\begin{equation*}
(x_0+x_2)^k \langle v' , \co(x_0) e^{x_2 D^*} v \rangle =(x_0+x_2)^k
\langle v', (Id_V \otimes e^{x_2 D^*})  \co(x_0 + x_2)v \rangle.
\end{equation*}

\noindent Since $(x_0+x_2)^k$, $\langle v' , \co(x_0)  e^{x_2 D^*} v
\rangle$ and $\langle v', (Id_V \otimes e^{x_2 D^*})  \co(x_0 + x_2)v
\rangle$ all involve no negative powers of $x_2$, we may multiply
both sides by $(x_0+x_2)^{-k}$ to get

\begin{equation*}
\langle v' , \co(x_0)  e^{x_2 D^*} v \rangle = \langle v', (Id_V
\otimes e^{x_2 D^*})  \co(x_0 + x_2)v \rangle
\end{equation*}

\noindent which is independent of $v'$ and $v$, or equivalently,

\begin{equation} \label{E:5221}
 (Id_V \otimes e^{-x_2 D^*}) \co(x_0)
=   \co(x_0 + x_2) e^{-x_2 D^*}.
\end{equation}

Finally, let $v' \in (V^{\otimes 3})'$, pick $k \in \Z_+$ so that
weak coassociativity holds for $v'(T \otimes Id_V)(Id_V \otimes T)$,
i.e. for any $v \in V$

\begin{multline} \label{E:5222}
(x_2-x_1)^k \langle v', (T \otimes Id_V)(Id_V \otimes T) (\co(x_2) \otimes Id_V)  \co(-x_1)v \rangle \\
=(x_2-x_1)^k \langle v', (T \otimes Id_V)(Id_V \otimes T) (Id_V
\otimes \co(-x_1))  \co(x_2 - x_1)v \rangle.
\end{multline}

\noindent Successively using skew-symmetry, (\ref{E:5222}),
(\ref{E:5221}), and skew-symmetry, we have

\begin{align*}
(x_2-&x_1)^k \langle v', (Id_V \otimes \co(x_2)) \co(x_1)v \rangle \\
&=(x_2-x_1)^k \langle v', (T \otimes Id_V) (Id_V \otimes T) (\co(x_2) \otimes Id_V) \co(-x_1) e^{x_1 D^*}v \rangle \\
&=(x_2-x_1)^k \langle v', (T \otimes Id_V)(Id_V \otimes T) (Id_V \otimes \co(-x_1))  \co(x_2 - x_1)e^{x_1 D^*}v \rangle \\
&=(x_2-x_1)^k \langle v', (T \otimes Id_V)(Id_V \otimes T) (Id_V \otimes \co(-x_1) e^{x_1 D^*} )  \co(x_2)v \rangle \\
&=(x_2-x_1)^k \langle v', (T \otimes Id_V)(Id_V \otimes \co(x_1))
\co(x_2)  v \rangle.
\end{align*}

\noindent This is a form of weak cocommutativity.
\end{proof}

The proposition allows us to consider weak coassociativity and
skew-symmetry as a replacement for the Jacobi identity in the
definition of graded vertex coalgebra.

\section{Comodules}
We now study comodules over vertex coalgebras.  Intuitively, we
desire to investigate a vector space under which left
comultiplication by a given vertex coalgebra satisfies all of the
axioms of vertex coalgebra that make sense.

\begin{definition}
A \emph{comodule} $M$ for a vertex coalgebra $(V,\co_V,c)$ is a vector space
equipped with a linear map

\begin{align*}
\co_M(x) : M &\mapsto (V \otimes M)[[x,x^{-1}]]
\end{align*}

\noindent such that for all $v \in V,$

 $\begin{array}{lc}
\text{1. Left counit:} & (c \otimes Id_M) \co_M(x)v=v \\
\text{2. Truncation:} & \co_M(x)v \in (V \otimes M)((x^{-1})) \\
\text{3. Jacobi Identity:}
 \end{array}
$

%
%

\begin{multline*} 
x_0^{-1}\delta \left(\frac{x_1-x_2}{x_0} \right) (Id_V \otimes
\co_M(x_2))  \co_M(x_1) -x_0^{-1}\delta \left(\frac{x_2-x_1}{-x_0}
\right)
(T \otimes Id_M)\\
 (Id_V \otimes \co_M(x_1))  \co_M(x_2)
=x_2^{-1}\delta \left(\frac{x_1-x_0}{x_2} \right) (\co_V(x_0)
\otimes Id_M)  \co_M(x_2).
\end{multline*}
\end{definition}


Moving to graded objects, we follow a precedent in vertex (operator)
algebras and allow comodules over graded vertex coalgebras to be
$\Q$-graded (cf. \cite{FHL}, although $\mathbb{C}$-grading is
sometimes used).

\begin{definition}
A $\Q$-graded vector space $M = \coprod_{j \in \Q} M_{(j)}$, such
that $M_{(j)}=0$ for $j<<0$, is considered a \emph{comodule over a
graded vertex coalgebra} $V$ if it is a comodule over $V$ as an
ungraded vertex coalgebra such that the map $\co_M (x) : M \mapsto
(V \otimes M)[[x,x^{-1}]]$ may be described by

\begin{equation*}
w \mapsto \co_M(x)w = \sum_{k\in \Z} \Delta_k(w) x^{-k-1}
\end{equation*}

\noindent where for $w \in M_{(j)}$, $\Delta_k(w) \in (V \otimes
M)_{(j+k+1)}$.  The graded vertex coalgebra will be said to be
\emph{finite dimensionally graded} if $\dim M_{(j)}< \infty$ for
each $j \in \Q$.
\end{definition}

Much of the work of studying graded vertex coalgebras extends to
comodules of graded vertex coalgebras. Referring to the proofs of
Propositions \ref{P:weakcom}, \ref{P:rr&com}, \ref{P:weakassoc} and
\ref{P:lr&ass} immediately give us the following propositions for
comodules.

\begin{proposition} \label{P:Mweakcom}
\textbf{\emph{(weak cocommutativity)}} Given a comodule $M$ over a
graded vertex coalgebra $V$, let $v' \in (V \otimes V \otimes M)'$.
There exists $k \in \Z_+$ such that

\begin{equation} \label{E:Mweakcom}
(x_1-x_2)^k \langle v', (Id_V \otimes \co_M(x_2))  \co_M(x_1)w -(T
\otimes Id_M) (Id_V \otimes \co_M(x_1))  \co_M(x_2)w \rangle=0
\end{equation}

\noindent for any $w \in M$.
\end{proposition}

In particular, we have the analog of Equation (\ref{E:5204}) for
comodules. For $v' \in V'$, $w' \in M'$, there exists some $K \in
\Z_+$

\begin{equation} \label{E:comodule_positive_powers}
\langle v' \otimes w', \co_M(x_0)w \rangle \in \C x_0^{-K}[[x_0]]
\end{equation}

\noindent for any choice of $w \in M$.

\begin{proposition}\label{P:Mrr&com} (a) \textbf{\emph{(right rationality)}}
Let a vector space $M$ satisfy all the axioms of being a comodule
over a graded vertex coalgebra $V$ except the Jacobi identity, but
also satisfy weak cocommutativity. Also let $v' \in (V \otimes V
\otimes M)'$ and $w \in M$. Then the formal series $\langle v',
(Id_V \otimes \co_M(x_2)) \co_M(x_1)w \rangle$ is in $\C[x_1
,x_2^{-1}][[x_1^{-1},x_2]]$ and in fact

\begin{equation*}
\langle v', (Id_V \otimes \co_M(x_2)) \co_M(x_1)w \rangle = \iota_{1
2} f(x_1,x_2)
\end{equation*}

\noindent where $f(x_1,x_2) \in \C[x_1,x_2]_S$ is uniquely
determined and is of the form

\begin{equation*}
f(x_1,x_2)= \frac{g(x_1,x_2)}{x_1^r x_2^s (x_1-x_2)^t}
\end{equation*}

\noindent for some $g(x_1,x_2) \in \C[x_1,x_2]$ and $r,s,t \in \Z$.

(b) \textbf{\emph{(cocommutativity)}} It is also the case that

\begin{equation} \label{E:Mcom}
\iota^{-1}_{1 2} \langle v', (Id_V \otimes \co_M(x_2)) \co_M(x_1)w
\rangle =\iota^{-1}_{2 1} \langle v', (T \otimes Id_M) (Id_V \otimes
\co_M(x_1)) \co_M(x_2)w \rangle.
\end{equation}
\end{proposition}

\begin{proposition} \label{P:Mweakassoc}
\textbf{\emph{(weak coassociativity)}} Given a comodule $M$ over a
graded vertex coalgebra $V$, let $v' \in (V \otimes V \otimes M)'$.
There exists $k \in \Z_+$ such that

\begin{equation} \label{E:Mweakassoc}
(x_0+x_2)^k \langle v', (\co_V(x_0) \otimes Id_M)  \co_M(x_2)w
-(Id_V \otimes \co_M(x_2))  \co_M(x_0 + x_2)w \rangle=0
\end{equation}

\noindent for any $w \in M$.
\end{proposition}

\begin{proposition}\label{P:Mlr&ass} (a) \textbf{\emph{(left rationality)}}
Let a vector space $M$ satisfy all the axioms of being a comodule
over a graded vertex coalgebra $V$ except the Jacobi identity, but
also satisfy weak coassociativity. Also let $v' \in (V \otimes V
\otimes M)'$ and $w \in M$. Then the formal series $\langle v',
(\co_V(x_0) \otimes Id_M) \co_M(x_2)w \rangle$ is in
$\C[x_0^{-1},x_2][[x_0,x_2^{-1}]]$ and in fact

\begin{equation*}
\langle v', (\co_V(x_0) \otimes Id_M) \co_M(x_2)w \rangle = \iota_{2
0} k(x_0,x_2)
\end{equation*}

\noindent where $k(x_0,x_2) \in \C[x_0,x_2]_S$ is uniquely
determined and is of the form

\begin{equation} \label{E:Mk_function}
k(x_0,x_2)= \frac{h(x_0,x_2)}{x_0^r x_2^s (x_0+x_2)^t}
\end{equation}

\noindent for some $h(x_0,x_2) \in \C[x_0,x_2]$ and $r,s,t \in \Z$.

(b) \textbf{\emph{(coassociativity)}} It is also the case that
\begin{equation} \label{E:Massoc}
\iota^{-1}_{1 2} \langle v', (Id_V \otimes \co_M(x_2)) \co_M(x_1)w
\rangle =\left. \left( \iota^{-1}_{2 0} \langle v', (\co_V(x_0)
\otimes Id_M) \co_M(x_2)w \rangle \right) \right|_{x_0=x_1-x_2}.
\end{equation}
\end{proposition}

The trivial modifications that we have made to right and left
rationality, cocommutativity and coassociativity allow us also to
claim Proposition \ref{P:RCAtoJI} in the graded vertex coalgebra
comodule context.

\begin{proposition}\label{P:MRCAtoJI}
The Jacobi identity may be replaced in the definition of a comodule
over graded vertex coalgebra by right and left rationality,
cocommutativity and coassociativity, i.e. the claims of Propositions
\ref{P:Mrr&com} and \ref{P:Mlr&ass}.
\end{proposition}

Unlike the case in graded vertex coalgebras, where the most natural
substitute for the Jacobi identity would seem to be weak
cocommutativity and the $D^*$-bracket, for comodules of graded vertex
coalgebras we are able to use only weak coassociativity to obtain
the Jacobi identity. (The $D^*$-bracket formula does not make sense in
the comodule context, and it is not known if cocommutativity itself
implies the Jacobi identity.)

\begin{proposition} \label{P:Mweakassoc_JI}
Let a vector space $M$ satisfy all the axioms of being a comodule
over a graded vertex coalgebra $V$ except the Jacobi identity, but
also satisfy weak coassociativity (\ref{E:Mweakassoc}). Then $M$
satisfies the Jacobi identity.
\end{proposition}

\begin{proof}  By Proposition \ref{P:Mlr&ass}, we have right
rationality and coassociativity.  We need only show weak
coassociativity of comodules implies right rationality and
cocommutativity, for then Proposition \ref{P:MRCAtoJI} implies the
Jacobi identity. Notice that right rationality is actually a direct
consequence of coassociativity (\ref{E:Massoc}) so our main task is
to show cocommutativity (\ref{E:Mcom}).

We start by showing that the $D^*$-derivative holds for comodule
comultiplication.  Let $v' \in V'$, $w' \in M'$ and pick $k \in
\Z_+$ so that weak coassociativity holds for $v' \otimes c \otimes
w'$. So for all $w \in M,$

\begin{multline*}
(x_0+x_2)^k \langle v' \otimes c \otimes w' , (\co_V(x_0) \otimes Id_M)  \co_M(x_2) w \rangle \\
=(x_0+x_2)^k \langle v' \otimes c \otimes w' , (Id_V \otimes
\co_M(x_2))  \co_M(x_0 + x_2)w \rangle.
\end{multline*}

\noindent Applying (\ref{E:cocreat3}) on the left-hand side and the
left counit property on the right, we get

\begin{equation*}
(x_0+x_2)^k \langle v' \otimes w' , (e^{x_0 D^*} \otimes Id_M)
\co_M(x_2) w \rangle =(x_0+x_2)^k \langle v' \otimes w' , \co_M(x_0
+ x_2) w \rangle.
\end{equation*}

\noindent
Recalling (\ref{E:comodule_positive_powers}), we may choose $k$ larger if
necessary so that $x^k \langle v' \otimes w'$, $\co_M(x)w \rangle \in \C [[x]]$.
For positive powers, $x_0 + x_2$ and $x_2 + x_0$ are equivalent, so we may write

\begin{equation*}
(x_0+x_2)^k \langle v' \otimes w' , (e^{x_0 D^*} \otimes Id_V) \co_M(x_2) w \rangle
=(x_0+x_2)^k \langle v' \otimes w' , \co_M(x_2 + x_0) w \rangle.
\end{equation*}

Finally, multiply both
sides by $(x_0+x_2)^{-k}$ and observe that independent of $v'$,
$w'$, $w$,

\begin{equation} \label{E:MDderiv}
(e^{x_0 D^*} \otimes Id_M)  \co_M(x_2) = \co_M(x_2 + x_0).
\end{equation}

Now by left rationality (a consequence of weak coassociativity), for
$v' \in (V \otimes V \otimes M)'$, $w \in M$ we have a rational function of
the form

\begin{equation*}
k(x_0,x_2)= \frac{h(x_0,x_2)}{x_0^r x_2^s (x_0+x_2)^t},
\end{equation*}

\noindent with $h(x_0,x_2) \in \C[x_0,x_2]$, $r,s,t \in \Z$ such that

\begin{equation} \label{E:314354}
\langle v', (\co_V(x_0) \otimes Id_M) \co_M(x_2)w \rangle = \iota_{2 0}
k(x_0,x_2).
\end{equation}

\noindent We know that by coassociativity, the left-hand side of the
cocommutativity equation (\ref{E:Mcom}) satisfies

\begin{align*}
\iota^{-1}_{1 2}\langle v', (Id_V \otimes \co_M(x_2)) \co_M(x_1)v \rangle
&= \left. \left( \iota^{-1}_{2 0} \langle v', (\co_V(x_0) \otimes Id_M) \co_M(x_2)w \rangle
\right) \right|_{x_0=x_1-x_2} \\
&= k(x_0,x_2)|_{x_0=x_1-x_2}.
\end{align*}

\noindent Now working from the right-hand side of cocommutativity
(\ref{E:Mcom}) and using coassociativity, (\ref{E:skew-symm}),(\ref{E:MDderiv}),
and (\ref{E:314354}) we also have

\begin{align*}
\iota^{-1}_{2 1}\langle v', (T \otimes Id_V)(&Id_V \otimes \co_M(x_1)) \co_M(x_2)  v \rangle \\
&= \left. \left( \iota^{-1}_{1 0} \langle v', (T\co_V(x_0) \otimes Id_M) \co_M(x_1)w \rangle
\right) \right|_{x_0=x_2-x_1} \\
&= \left. \left( \iota^{-1}_{1 0} \langle v', (\co_V(-x_0)e^{x_0 D^*}
\otimes Id_M) \co_M(x_1)w \rangle
\right) \right|_{x_0=x_2-x_1} \\
&= \left. \left( \iota^{-1}_{1 0} \langle v', (\co_V(-x_0) \otimes
Id_M) \co_M(x_1+x_0)w \rangle
\right) \right|_{x_0=x_2-x_1} \\
&= k(-x_0,x_1+x_0)|_{x_0=x_2-x_1}.
\end{align*}

\noindent But as rational functions,

\begin{equation*}
k(x_0,x_2)|_{x_0=x_1-x_2} = \frac{h(x_1-x_2,x_2)}{(x_1-x_2)^r x_2^s
x_1^t} = k(-x_0,x_1+x_0)|_{x_0=x_2-x_1}.
\end{equation*}

\noindent This completes the proof of cocommutativity.
\end{proof}

Hence, comodules of graded vertex coalgebras may now be axiomatized
as follows.

\begin{definition}
A \emph{comodule} over a graded vertex coalgebra $V$ may be
described as a collection of data

\noindent 1. a $\Q$-graded vector space $V=\coprod_{k \in \Q}
V_{(k)}$, with $V_{(k)}=0$ for $k <<0$ \\
2. a linear map

\begin{align*}
\co_M (x) : M &\to (V \otimes M)((x^{-1})) \\
v &\mapsto \co_M(x)v = \sum_{k\in \Z} \Delta_k(v) x^{-k-1},
\end{align*}

such that for $v \in V_{(r)}$, $\Delta_k(v) \in (V \otimes M)_{(r+k+1)}$, \\

\noindent
satisfying the following axioms for $v' \in (V \otimes V \otimes M)'$ and $w \in M$: \\
1. $(c \otimes Id_M) \co_M(x)=Id_M$ \\
2. There exists $k \in \Z_+$ depending on $v'$ but independent of
$w$ such that

\begin{equation*}
(x_0+x_2)^k \langle v', (\co_V(x_0) \otimes Id_M)  \co_M(x_2)w -
(Id_V \otimes \co_M(x_2))  \co_M(x_0 + x_2)w \rangle=0.
\end{equation*}
\end{definition}

\section{Vertex operator coalgebras and their comodules}

Starting with \cite{Hub1} and continuing in \cite{Hub2}, \cite{Hub3}
and \cite{Hub4}, a particular type of graded vertex algebra, called
a vertex operator coalgebra, has been studied which possesses a
conformal structure. We note that the cocommutativity,
coassociativity and comodule results for graded vertex coalgebras
extend automatically to vertex operator coalgebras.

A vertex operator coalgebra may be thought of as a vertex algebra
with a natural action of the Virasoro algebra.  Natural in this
context means that the grading of the vertex operator algebra is the
one determined by the generator $L(0)$ of the Virasoro algebra and
$L(1)$ acts as the operator $D^*$.

%
%
%
%
%
%
%
%
%
%
%
%
%
%
%
Consequently, Propositions \ref{P:skew-symm}, \ref{P:Dbracket},
\ref{P:weakcom}, \ref{P:rr&com}, \ref{P:weakassoc}, \ref{P:lr&ass},
\ref{P:weakcom_JI} and \ref{P:weakassoc_JI} all hold for vertex
operator algebras.  Not only do these propositions give insight
into the structure of vertex operator coalgebras, they also justify the following
equivalent definition.

\begin{definition}
A \emph{vertex operator coalgebra of rank $d \in \C$} has the
following data:

\begin{itemize}

\item $V = \coprod_{k = -N}^{\infty} V_{(k)}$, with $\dim V_{(k)} < \infty$, $N \in \N$

\item
a linear map $c : V \mapsto \C$,

\item
a linear map $\rho : V \mapsto \C$,

\item
a linear map
\begin{align*}
\co (x) : V &\mapsto (V \otimes V)((x^{-1})) \\
v &\mapsto \co(x)v = \sum_{k\in \Z} \Delta_k(v) x^{-k-1}
\end{align*}
\end{itemize}

satisfying (for $v \in V$):

\begin{tabbing}
1. Left counit: \ \ \ \ \ \ \
\= $(c \otimes Id_V) \co(x)v=v$ \\
2. Cocreation:
\> $(Id_V \otimes c) \co(x)v \in V[[x]]$ \\
\> \hspace{-.5 in} $\lim_{x \to 0} (Id_V \otimes c) \co(x)v=v$ \\
\end{tabbing}

\noindent 3. Weak cocommutativity: for $v' \in (V^{\otimes 3})'$,
there exists $k \in \Z_+$ (depending on $v'$ but not on $v$) such
that
$$
(x_1-x_2)^k \langle v',  (Id_V \otimes \co(x_2))  \co(x_1)v - (T \otimes Id_V) (Id_V \otimes \co(x_1))
\co(x_2)v \rangle=0
$$
4. Virasoro algebra:
given
$(\rho \otimes Id_V) \co(x) = \sum_{k \in \Z} L(k) x^{k-2}$, for all $j, k \in \Z$
\begin{equation*}
[L(j),L(k)]=(j-k)L(j+k)+\frac{1}{12}(j^3-j)\delta_{j,-k}d,
\end{equation*}
5. Grading:  \ \
 if $v \in V_{(k)}$, $L(0)v= kv$ \\
6. $L(1)$-Derivative and bracket:

\begin{equation*}
\frac{d}{dx} \co(x) = (L(1) \otimes Id_V) \co(x)
=\co(x)L(1) -(Id_V \otimes L(1)) \co(x).
\end{equation*}
\end{definition}

In practice weak cocommutativity is easier to verify than the Jacobi
identity, making this description preferable in contexts where verification is necessary.

A comodule over a vertex operator coalgebra has a similar alternate definition as a result of
Propositions \ref{P:Mweakcom}, \ref{P:Mrr&com}, \ref{P:Mweakassoc},
\ref{P:Mlr&ass}, and \ref{P:Mweakassoc_JI}
%
%
%
%
%
%
%
%
%
%
%

\begin{definition}
A \emph{vertex operator coalgebra comodule of rank $d \in \C$} over
$V$ has the following data:

\begin{itemize}

\item $M = \coprod_{k \in \Q} M_{(k)}$, with $\dim M_{(k)} <
\infty$ for all $k$, $M_{(k)}=0$ for $k << 0$,

\item
a linear map
\begin{align*}
\co_M (x) : M &\mapsto (V \otimes M)((x^{-1})) \\
w &\mapsto \co_M(x)w = \sum_{k\in \Z} \Delta_k(w) x^{-k-1}
\end{align*}
\end{itemize}

satisfying (for $w \in M$):

\begin{tabbing}
1. Left counit: \ \ \ \ \ \ \
\= $(c \otimes Id_V) \co_M(x)w=w$
\end{tabbing}

\noindent 2. Weak coassociativity: for $v' \in (V \otimes V \otimes
M )'$, there exists $k \in \Z_+$ depending on $v'$ (but not on $w$)
such that
$$
(x_0+x_2)^k \langle v', (\co_V(x_0) \otimes Id_M)  \co_M(x_2)w
- (Id_V \otimes \co_M(x_2))  \co_M(x_0 + x_2) w
\rangle=0
$$
\noindent 3. Virasoro algebra: given $(\rho \otimes Id_M) \co_M(x) =
\sum_{k \in \Z} L(k) x^{k-2}$, for all $j, k \in \Z$
\begin{equation*}
[L(j),L(k)]=(j-k)L(j+k)+\frac{1}{12}(j^3-j)\delta_{j,-k}d,
\end{equation*}
4. Grading:  \ \ \ \ \ \ \ \ \ \ \
 if $w \in M_{(k)}$, $L(0)w= kw$ \\
5. $L(1)$-Derivative and bracket:

\begin{align*}
\frac{d}{dx} \co_M(x) &= (L(1) \otimes Id_M) \co_M(x) \\
&=\co_M(x)L(1) -(Id_V \otimes L(1)) \co_M(x).
\end{align*}
\end{definition}

\section{The categories of finite dimensionally graded vertex coalgebras and finite dimensionally
graded vertex algebras}

An additional benefit of studying finite dimensionally graded vertex
coalgebras is that grading allows us to pass naturally between a
vector space and it's graded dual, then back again while the finite
dimensionality prevents enlarging our vector space. First, we
restrict our attention to graded vertex coalgebras as described
above and to the category graded vertex algebras, by which we mean a
vertex algebra $(V,Y,\mathbf{1})$ such that $V = \coprod_{k \in \Z}
V_{(k)}$ with all $V_{(k)}$ trivial for $k << 0$, and for each $u
\in V_{(j)}$, $v \in V_{(K)}$, we have $u_n v \in V_{(j+k-n-1)}$.

Here we mainly follow the approach in \cite{Hub4}.

\begin{theorem} \label{T:bij_VA_VC}
Let $V = \coprod_{k \in \Z} V_{(k)}$ be a $\Z$-graded vector space.
Choose a distinguished vectors $\mathbf{1} \in V_{(0)}$ and a linear
map

\begin{equation*}
Y(\cdot, x): V \to \text{End }V [[x,x^{-1}]].
\end{equation*}

Additionally, define $c \in V''_{(0)}$ to be the double dual of
$\mathbf{1}$, and a linear operator $\co (x): V' \to (V' \otimes
V')[[x,x^{-1}]]$ defined by

\begin{equation*}
\langle \co(x)u', v \otimes w \rangle = \langle u', Y(v, x)w \rangle
\end{equation*}
for all $u' \in V'$, $v,w \in V$. Then $(V, Y, \mathbf{1})$ is a
graded vertex algebra if and only if $(V', \co, c)$ is a graded
vertex coalgebra.
\end{theorem}

The proof follows that of Theorem 2.9 in \cite{Hub4} quite neatly,
but with two interesting features.  First, the Virasoro bracket
structure is not intermingled with the proof of any of the other
axioms and thus may be entirely removed.  Second, the gradation, or
specifically the weight condition $u_n v \in V_{(j+k-n-1)}$ for $u
\in V_{(j)}$, $v \in V_{(K)}$, is essential to proving the
truncation condition.  The proof of Theorem \ref{T:bij_VA_VC}
depends on having \emph{graded} vertex algebras and coalgebras.

As in \cite{Hub4}, Theorem \ref{T:bij_VA_VC} gives rise to a
contravarient functor from the category of graded vertex algebras to
the category of graded vertex coalgebras, and a similar functor goes
in the opposite direction.  When $V$ is a finite dimensionally
graded vertex algebra, $V$ and $V''$ are canonically isomorphic, and
these functors demonstrate an equivalence of categories.

\begin{corollary}
The categories of finite dimensionally graded vertex algebras and
finite dimensionally graded vertex coalgebras are equivalent.
\end{corollary}

\section{Motivation from the universal enveloping algebra of a finite dimensional Lie algebra} \label{S:associative}

Vertex algebras are distinctive in that they exhibit properties of
associative algebras and of Lie algebras.  In trying to determine
the correct ``bialgebraic" structures appearing in the vertex
algebra world, it is productive to consider the most widely known
example of bialgebras - the universal enveloping algebra of a Lie
algebra. Given a finite dimensional abelian Lie algebra
$\mathcal{L}$, we construct the universal enveloping algebra of
$\mathcal{L}$ by first considering the tensor algebra of
$\mathcal{L}$ as a vector space
$$
T(\mathcal{L})= \coprod_{i \in \N} \mathcal{L}^{\otimes i}
$$
where $\mathcal{L}^{\otimes i}= \underbrace{\mathcal{L} \otimes
\cdots \otimes \mathcal{L}}_i$ (and $\mathcal{L}^{\otimes 0}=\C$).
(We could use a Heisenberg algebra, or any Lie algebra where the
center contains the commutator, but an abelian Lie algebra will do
for the purpose of the current analogy.) Let $J$ be the ideal
generated by the tensors
$$
x \otimes y -y \otimes x - [x,y]
$$
for $x,y\in \mathcal{L}$.  The universal enveloping algebra is the
associative, coassociative bialgebra
$$
U(\mathcal{L})=T(\mathcal{L})/J
$$
where multiplication is given by
$$
m(u,v)=u \otimes v
$$
for $u,v \in U(\mathcal{L})$.

Now any coalgebra structure on a vectors space $V$ naturally gives
rise to an algebra structure on the dual space $V^*$, i.e. for
$u',v'\in V^*$, $w\in V$, $m(u',v')(w)=(u' \otimes v') (\Delta w)$.
But an algebra structure on a vector space $V$ only gives a
coalgebra structure on a subspace of the dual,
$$
V^{\circ}=\{f \in V^* | f(I)=0 \text{ for some ideal } I \text{ of } V
\text{ such that dim } V/I < \infty \},
$$
the so-called `finite dual' (cf. \cite{M}).  The comultiplication
$\Delta^*:V^{\circ} \to V^{\circ} \otimes V^{\circ} $ is determined by
$$
\Delta^*(u')(v \otimes w)=u'(m(v,w))
$$
for $u'\in V^{\circ}$, $v,w\in V$, where $m$ is the multiplication on $V$.

In order for $U(\mathcal{L})$ to be a bialgebra however, we must have a comultiplication on the vector space
itself rather than on the finite dual.  We will construct the comultiplicative structure on $U(\mathcal{L})$
as follows.  By the
Poincar\'e-Birkhoff-Witt theorem, a basis $e_1, \ldots, e_n$ of $\mathcal{L}$ determines
a basis from $U(\mathcal{L})$, which is
$$
\{ e_{i_1} \otimes \cdots \otimes e_{i_k} | k \in \N, i_1 \leq i_2 \leq \cdots \leq i_k \leq n \}
$$
We then have a linear map
$$\Phi: U(\mathcal{L}) \to U(\mathcal{L})^{\circ}$$
determined by mapping each basis element $e_{i_1} \otimes \cdots
\otimes e_{i_k}$ to the dual element $f \in U(\mathcal{L})^*$ which
maps $e_{i_1} \otimes \cdots \otimes e_{i_k}$ to one and all other
basis elements to 0.  Certainly $f \in U(\mathcal{L})^{\circ}$
since allowing
$$
I=\text{span}\{ e_{j_1} \otimes \cdots \otimes e_{j_\ell} | \ell >
k, j_1 \leq j_2 \leq \cdots \leq j_{\ell} \leq n \},
$$
we have $f(I)=0$ and dim$\left( U(\mathcal{L})/I\right)< \infty$.
The map $\Phi$ is definitely injective, and surjectivity may be seen
as follows.  Any $f \in U(\mathcal{L})^*$ has a set of basis
elements $e_{i_1} \otimes \cdots \otimes e_{i_k}$ for which
$f(e_{i_1} \otimes \cdots \otimes e_{i_k})\neq 0$. If the set is
finite, $f \in \Phi[U(\mathcal{L})].$ If the set is infinite, it is
impossible to pick an ideal $I$ of $U(\mathcal{L})$ such that
$f(I)=0$ and dim $U(\mathcal{L})/I < \infty$. Hence, $\Phi$ is a
vector space isomorphism.

Using
the isomorphism $\Phi$, we induce a comultiplication on
$U(\mathcal{L})$ from the comultiplication $\Delta^*$ on $U(\mathcal{L})^{\circ}$,
$$
\Delta=(\Phi \otimes \Phi)^{-1} \circ \Delta^* \circ
\Phi.
$$
We will see that this comultiplication is the `classical' (and
unique, cf. \cite{U}) comultiplication $\Delta$ on $U(\mathcal{L})$.
We need only examine $\Delta$ on $\mathcal{L} \subset
U(\mathcal{L})$, since the comultiplication is completely determined
by these elements using the properties:
$$
\Delta(1)=1 \otimes 1,
$$
$$
\Delta(vw)=(m \otimes m) \circ (Id_{U(\mathcal{L})} \otimes T
\otimes Id_{U(\mathcal{L})}) (\Delta(v) \otimes \Delta(w)),
$$
for any $v, w \in U(\mathcal{L})$.  First, consider a basis element
$u \in \mathcal{L} \subset U(\mathcal{L})$, and apply $\Delta^*
\circ \Phi (u)$ to two arbitrary basis elements $v, w \in
U(\mathcal{L})$. Denote by $u' \in U(\mathcal{L})^{\circ}$ the
element that is 1 on $u$ and zero on every other basis element.
\begin{align*}
\Delta^* \circ \Phi (u)(v \otimes w)
&= \Delta^* (u')(v \otimes w)\\
&= u'(v w) \\
&= \left\{
\begin{array}{cc}
1,&\mbox{ if } u=vw\\
0, &\mbox{ otherwise }
\end{array}\right.
\end{align*}
Thus $\Delta^* \circ \Phi (u)$ is the element of
$U(\mathcal{L})^{\circ} \otimes U(\mathcal{L})^{\circ}$ that takes a
value of 1 on the basis elements $u \otimes 1$ and $1 \otimes u$,
but takes a value of zero on every other basis element.  In other
words, given $u'$ the dual of $u$ as above and $1'$ the dual of 1,
$$
\Delta^* \circ \Phi (u)= u' \otimes 1' + 1' \otimes u'.
$$
Therefore
$$
\Delta(u)= u \otimes 1 + 1 \otimes u.
$$
This is the classical description of the comultiplication of a
universal enveloping algebra, but we have derived it by examining
the comultiplication induced on the finite dual by the original
\emph{multiplication}.

With this as motivation, we now turn our attention the vertex
algebras and repeat the process to find a `correct' comultiplicative
structure.  Rather than posing this result in terms of dual vector
spaces or finite duals, these results may be posed in terms of a
nondegenerate bilinear form $( \cdot, \cdot ):V \otimes V \to \C$.
Given any vector space $V$, a bilinear form on $V$ determines a
linear map from $V$ to $V^*$, and vice versa. A nondegenerate
bilinear form corresponds to an injective map $\Phi: V \to V^*$.
Thus, our isomorphism $\Phi:U(\mathcal{L}) \to
U(\mathcal{L})^{\circ} \subset U(\mathcal{L})^*$ determines a
nondegenerate bilinear form on $U(\mathcal{L})$ and we may define
$\Delta$ by
$$
( \Delta(u), v \otimes w ) = ( u, m(v,w))
$$
for all $u,v,w \in U(\mathcal{L})$.

We now employ the same technique on a graded vertex algebra $(V, Y,
\mathbf{1})$. For each $k \in \Z$, choose a basis $\{e_{k,i} \}_{i
\in I_k}$ for $V_{(k)}$.  Then $\{e_{k,i} | k \in \Z, i\in I_k \}$
is a basis for $V$ and determines a bilinear form $(\cdot,\cdot):V
\otimes V \to \C$ on $V$ by the rule
\begin{equation} \label{E:bilin}
(e_{k,i},e_{\ell,j})=\delta_{k, \ell} \delta_{i, j}.
\end{equation}
In other words, the product of any two basis elements is 1 if they
are the same and 0 if they are different.  Hence the form is
symmetric and nondegenerate by construction.  It is also graded,
i.e.,
$$(V_{(k)},V_{(\ell)})=0
$$
for $k\neq\ell$.

\begin{theorem} \label{T:gva_gvc_struc}
Given a graded vertex algebra $(V,Y, \mathbf{1})$ equipped with a
graded, nondegenerate bilinear form, $V$ also carries the structure
of a vertex coalgebra with linear operators
$$
c: V \to \C, \ v \mapsto (v,\mathbf{1})
$$
and
$$
\co(x):V \to (V \otimes V)[[x,x^{-1}]]
$$
defined by $(\co(x)u,v \otimes w)= (u, Y(v,x)w)$ for $u,v,w \in V$.
\end{theorem}

The proof follows that of Theorem 3.1 in \cite{Hub3} (which was
stated in terms of vertex operator algebras and coalgebras) with
weakened requirements. Instead of requiring a form that is invariant
or Virasoro preserving, a graded form is all that is used in any
portion of the proof besides the Virasoro bracket.

By the form constructed directly preceding Theorem
\ref{T:gva_gvc_struc}, we conclude:

\begin{corollary} \label{C:main_cor}
Every graded vertex algebra has a graded vertex coalgebra structure via
the nondegenerate form determined in Equation (\ref{E:bilin}).
\end{corollary}

In the final section, we will extend the parallel with (associative)
enveloping algebras by examining enveloping vertex algebras.

\section{A graded vertex coalgebra structure on the universal enveloping algebra of vertex Lie
algebras} \label{S:Enveloping_Vertex_Algebras}

In \cite{P}, Primc defines the notion of vertex Lie algebra, which
contains the ``positive" information of a vertex algebra.  He then
constructs an enveloping vertex algebra for a vertex Lie algebra and
goes on to show that the enveloping vertex algebra has the same
universal property as that of an associative enveloping algebra.
Namely, consider a (vertex) Lie Algebra $\mathcal{L}$, a (vertex)
algebra $A$, and a (vertex) Lie algebra homomorphism
$\phi:\mathcal{L} \to A$, viewing $A$ as a (vertex) Lie algebra.
There is a unique embedding $\iota:\mathcal{L} \to U(\mathcal{L})$
from $\mathcal{L}$ to its enveloping (vertex) algebra and a unique
(vertex) algebra homomorphism $\psi: U(\mathcal{L}) \to A$, such
that $\psi \circ \iota = \phi$ (cf. \cite{Hum}).

The classical Lie algebra's enveloping algebra has the structure not
only of a bialgebra, but a Hopf algebra.  In the vertex Lie algebra
setting, a general vertex bialgebra or vertex Hopf algebra structure
has yet to be defined. However, we will show that given a graded
vertex Lie algebra, Primc's enveloping vertex algebra is naturally a
grade vertex coalgebra. This is a step toward giving the enveloping
vertex algebra a general vertex bialgebra structure.  (In \cite{L2},
Li defines a structure with a vertex algebra and a coassociative
coalgebra structure, but points out there are more general
constructions.)

We begin by rigorously defining vertex Lie algebra, adding a grading
and constructing the (now graded) enveloping vertex algebra. Our
definition will make use of the symbol $\simeq$ which indicates
equality of the principal parts of two formal Laurent series. For
example, given two Laurent series $f(x_1,x_2)$ and $g(x_1,x_2)$ in
formal variables $x_1$ and $x_2$, $f \simeq g$ means that for all
$j,k \in \Z_+$ the coefficient of $x_1^{-j}x_2^{-k}$ in $f$ is equal
to the same coefficient in $g$.  Note below that in addition to
adding grading, we are restricting Primc's work over vertex
superalgebras to the vertex algebra setting.

\begin{definition}\label{D:voa}
A \emph{vertex Lie algebra} is a vector space $V$ equipped with a
linear operator $D:V \to V$ (called the derivation) and a linear map

\begin{align*}
Y (\cdot,x): V &\to (\text{End } V) [[x^{-1}]] \\
v &\mapsto Y(v,x) = \sum_{k\in \N} v_k x^{-k-1} \ \ \ \ \text{
(where } v_n \in \text{End }V),
\end{align*}

\noindent satisfying the following axioms for $u,v \in V$:

1. Truncation: $u_k v = 0$ for $k$ sufficiently large.

2. Half Jacobi identity:

\begin{multline} \label{E:Jac2}
x_0^{-1}\delta \left(\frac{x_1-x_2}{x_0} \right) Y(u,x_1)Y(v,x_2)
-x_0^{-1}\delta \left(\frac{x_2-x_1}{-x_0} \right)
Y(v,x_2)Y(u,x_1) \\
\simeq x_2^{-1}\delta \left(\frac{x_1-x_0}{x_2} \right)
Y(Y(u,x_0)v,x_2).
\end{multline}

3. Half skew-symmetry:  $Y(u,x)v \simeq e^{xD} Y(v,-x)u$.

4. D-bracket: $[D,Y(u,x)]=Y(Du,x)=\frac{d}{dx}Y(u,x)$.

\end{definition}

Certainly any vertex algebra is a vertex Lie algebra with
$D(v)=v_{-2}\mathbf{1}$. A vertex Lie algebra may also carry a
natural $\Z$-grading.

\begin{definition}
Let $V = \coprod_{k \in \N} V_{(k)}$ be a $\Z$-graded vector space..
A vertex Lie algebra structure on $V$ will be said to be
\emph{graded} if for each $r,s, t \in \N$, $u \in V_{(r)}$, and $v
\in V_{(s)}$,
$$
u_t v \in V_{(r+s-t-1)},
$$
and
$$
Du \in V_{(r+1)}.
$$
\end{definition}
This definition agrees with the notion of grading for vertex
(operator) algebras.  Now consider the affinization of a vertex Lie
algebra $V$ with derivation $D$, written $V \otimes \C[q,q^{-1}]$.
For graded vertex Lie algebra, $V \otimes \C[q,q^{-1}]$ carries a
natural $\Z$-grading by letting $q^{k}$ have weight $-k-1$.
Notationally, given $v\in V$ and $n \in \Z$, we use the abbreviation
$$
v_n=v \otimes q^n.
$$
Primc shows that the quotient space
$$
\mathcal{L}(V)=\left(V \otimes \C[q,q^{-1}] \right) / \{(Dv)_n+ n \
v_{n-1} | v \in V, \ n \in \Z \}
$$
is a Lie algebra under the bracket
$$
[u_m,v_n]= \sum_{i\geq0} {n \choose i} (u_i v)_{m+n-i},
$$
for $u,v \in V$ and $m,n \in \Z$, where the images of $u_m, v_n \in
V \otimes \C[q,q^{-1}]$ in $\mathcal{L}(V)$ are again denoted $u_m$
and $v_n$, respectively. (Employing the notation of \cite{P}, $(u_i
v)_{m+n-i}$ means the image $u_iv \otimes q^{m+n-i}$.)

Since $(Dv)_n$ has the same weight as $v_{n-1}$ for any homogeneous
$v\in V$, $\mathcal{L}(V)$ is $\Z$-graded.  Further,
$D:\mathcal{L}(V) \to \mathcal{L}(V)$ given by $D(v_n)=(Dv)_n$ is a
grading preserving derivation of the Lie algebra. Finally, given
$k,\ell,m,n \in \Z$, $u \in V_{(k)}$ and $v \in V_{(\ell)}$, $(u_i
v)_{m+n-i}$ has weight $k+\ell-m-n-2$ for any $i$, which is the
weight that $[u_m,v_n]$ must have in order to make the quotient
$\mathcal{L}(V)$ a $\Z$-graded Lie algebra with a grading preserving
derivation.

Moving forward, $\mathcal{L}(V)$ splits into two $D$-invariant Lie
subalgebras,
\begin{align*}
\mathcal{L}_-(V) &= \text{span} \{v_n|v \in V, n < 0\},\\
\mathcal{L}_+(V) &= \text{span} \{v_n|v \in V, n \geq 0\}.
\end{align*}
The enveloping vertex algebra we seek is a generalized Verma module
induced from a trivial $\mathcal{L}_+(V)$-module $\C$:
$$
\mathcal{V}(V)= U(\mathcal{L}(V)) \otimes_{U(\mathcal{L}_+(V))} \C,
$$
where $U$ gives the enveloping algebra of a given Lie algebra. Primc
established that $\mathcal{V}(V)$ is a vertex algebra, with the
derivation $D$ on $\mathcal{L}(V)$ extending naturally to
$\mathcal{V}(V)$ and with vacuum vector, $\mathbf{1},$ the tensor of
$1 \in U(\mathcal{L}(V))$ and $1 \in \C$.

Again, if $V$ is $\Z$-graded then the grading is preserved in $
\mathcal{V}(V)$ and we write
$$
\mathcal{V}(V)= \coprod_{k\in \Z} \mathcal{V}(V)_{(k)}.
$$
Exploiting the natural vector space isomorphism
$$
\mathcal{V}(V) \cong U(\mathcal{L}_-(V)),
$$
the elements of $\mathcal{V}(V)$ may be written
as sums of homogeneous elements of the form
$$
\hat{v}=v^{(1)}_{n_1} v^{(2)}_{n_2}\cdots v^{(t)}_{n_t} \in
\mathcal{V}(V)_{(p_1+\cdots+p_t-n_1-\cdots-n_t-t)}
$$
where $t\in \N$; $n_1,n_2,\ldots,n_t \in \Z_{-}$; $v^{(1)}\in
V_{(p_1)}, v^{(2)}\in V_{(p_2)}, \ldots, v^{(t)}\in V_{(p_t)}$ for
some $p_1, p_2, \ldots, p_t \in \N$; and tensor products are omitted
from the notation.  The derivation and the product operation both
naturally preserve the grading, 
as they are both induced from $\mathcal{L}(V)$. (This motivates the
above choice of grading: wt $v_n =$ wt $v -n-1$.) Thus
$\mathcal{V}(V)$ is a graded vertex algebra and by Corollary
\ref{C:main_cor}, $\mathcal{V}(V)$ also carries the structure of a
coalgebra.

The map $\iota_V: V \to \mathcal{L}_-(V)$ given by
$\iota_V(v)=v_{-1}$ is a vector space isomorphism by Theorem 4.6 in
\cite{P}.  Thus, if $V$ is finite dimensionally graded, as is the
case for vertex operator algebras, $\mathcal{L}_-(V)$ is finite
dimensionally graded and as a consequence $\mathcal{V}(V)$ is a
finite dimensionally graded vertex algebra and coalgebra.   By the
discussion at the beginning of Section
\ref{S:Enveloping_Vertex_Algebras}, it is reasonable to conclude
that this is the `appropriate' multiplicative structure for an
enveloping vertex bialgebra, and by the discussion in Section
\ref{S:associative}, it is reasonable to conclude that this is the
appropriate comultiplicative structure.


\end{document}